\documentclass{amsproc}

\usepackage{a4wide}
\usepackage{enumerate}
\usepackage{stmaryrd}
\usepackage{bbm}

\font\cmssl=cmss10 at 10 pt  

\newtheorem{theorem}{Theorem}[section]
\newtheorem{lemma}[theorem]{Lemma}
\newtheorem{Prop}[theorem]{Proposition}
\newtheorem{Cor}[theorem]{Corollary}
\theoremstyle{definition}
\newtheorem{definition}[theorem]{Definition}

\theoremstyle{remark}
\newtheorem{remark}[theorem]{Remark}

\numberwithin{equation}{section}

\newcommand{\hook}{\lrcorner \,}

\newcommand{\ve}{\varepsilon}
\newcommand{\Stab}{{\rm Stab}}

\newcommand{\aeh}{almost $\varepsilon$-Hermitian~}

\newcommand{\pJ}{J}  
\newcommand{\id}   {{\mathbbm{1}}}   
\newcommand{\eC}{{\mathbb C_\varepsilon}}
\newcommand{\hati}{{\rm i}_\varepsilon}
\newcommand{\nJ}[2]{(\n_{#1}\pJ)\,#2}
\newcommand{\no}[3]{(\n_{#1} \omega)\,(#2,#3)}
\newcommand{\metric}[2]{g( #1,#2 )}

\newcommand{\g}{\mathfrak g}
\newcommand{\h}{\mathfrak h}

\newcommand{\bt}{\begin{theorem}\ \ }  
\newcommand{\et}{\end{theorem}}  
\newcommand{\bp}{\begin{Prop}\ \ }  
\newcommand{\ep}{\end{Prop}}  
\newcommand{\bc}{\begin{Cor}\ \ }  
\newcommand{\ec}{\end{Cor}}  
\newcommand{\bl}{\begin{lemma}\ \ }  
\newcommand{\el}{\end{lemma}}  
\newcommand{\bd}{\begin{definition}\ \ }  
\newcommand{\ed}{\end{definition}}  
\newcommand{\pf}{\begin{proof}}  
\newcommand{\epf}{\end{proof}}  
\newcommand{\br}{\begin{remark}\ \ }
\newcommand{\er}{\end{remark}}
\newcommand{\brsn}{\begin{remarks*}\ \ }
\newcommand{\ersn}{\end{remarks*}}

\newcommand{\GL}{{\rm GL}}
\newcommand{\SL}{{\rm SL}}
\newcommand{\SO}{{\rm SO}}

\newcommand{\U}{{\rm U}}
\newcommand{\SU}{{\rm SU}}
\newcommand{\Mat}{{\rm Mat}}

\newcommand{\be}{\begin{equation}}
\newcommand{\ee}{\end{equation}}

\newcommand{\arr}{\begin{array}{rlll}}
\newcommand{\ea}{\end{array}}
\newcommand{\bea}{\begin{eqnarray}}
\newcommand{\eea}{\end{eqnarray}}
\newcommand{\bean}{\begin{eqnarray*}}
\newcommand{\eean}{\end{eqnarray*}}

\newcommand{\op}{\oplus}
\newcommand{\ot}{\otimes}
\newcommand{\ra}{\rightarrow}

\newcommand{\n}{\nabla}
\newcommand{\bn}{\bar \nabla} 

\newcommand{\bR}{\mathbb{R}}

\renewcommand{\o}{\omega}

\begin{document}

\title{Nearly pseudo-K\"ahler and nearly para-K\"ahler six-manifolds}

%    Information for first author
\author{Lars Sch\"afer}
\address{Lars Sch\"afer, Institut Differentialgeometrie, Leibniz Universit\"at
  Hannover, Welfengarten 1, D-30167 Hannover, Germany}
\email{schaefer@math.uni-hannover.de}

%    Information for second author
\author{Fabian Schulte-Hengesbach}
\address{Fabian Schulte-Hengesbach, Department Mathematik, Universit\"at Hamburg, Bundesstra{\ss}e 55, 
D-20146 Hamburg, Germany}
\email{schulte-hengesbach@math.uni-hamburg.de}

\date{\today}

\begin{abstract}
The subject of this paper is six-dimensional nearly (para-)K\"ahler geometry with pseudo-Riemannian metrics. Firstly, we derive the analogue of the well-known exterior differential system characterising a nearly K\"ahler manifold and prove applications to the automorphism group of a nearly (para-)K\"ahler structure. Secondly, we prove existence and uniqueness results for left-invariant nearly (para-)K\"ahler structures on Lie groups $G \times G$ where $G$ is three-dimensional and simple.
\end{abstract}

\maketitle
\setcounter{tocdepth}{2}

\section{Introduction}
The notion of a nearly K\"ahler manifold was introduced and studied in a series
of papers by A.\ Gray in the seventies in the context of weak holonomy. In the
last two decades, six-dimensional nearly K\"ahler manifolds turned out to be
of interest in a multitude of different areas including $\SU(3)$-geometries, stable forms, geometries with torsion, existence of Killing spinors, (weak) holonomy, supersymmetric models and compactifications of string theories. For a survey explaining the relations between most of these areas we refer to \cite{Ag}.

One observes that most of the literature on nearly K\"ahler geometry deals with Riemannian signature. To our best knowledge the paper on 3-symmetric spaces \cite{G2} is the only article by Gray considering also indefinite nearly K\"ahler metrics. 
Killing spinors on pseudo-Riemannian manifolds were studied in \cite{Ka} where nearly pseudo-K\"ahler and nearly para-K\"ahler manifolds appear in a natural way. 
The subject of nearly para-K\"ahler manifolds was further developed in \cite{IZ}. The prefix ``para'' roughly means that the anti-involutive complex structure is replaced by an involutive para-complex structure. We refer to section  \ref{sect_alm} for details on para-complex geometry.

Motivated by a class of solutions of the
topological-antitopological fusion equations on the tangent bundle
\cite{S1,S2} and the
similarity to special K\"ahler geometry, we became
interested in Levi-Civita flat nearly K\"ahler and Levi-Civita flat nearly para-K\"ahler manifolds. A classification of these manifolds in a constructive manner has been established in \cite{CS,CS2}. From these results 
it follows, that non-K\"ahlerian examples only exist in pseudo-Riemannian 
geometry. 
In other words, nearly K\"ahler geometry in the
pseudo-Riemannian world can be very different from the better-understood Riemannian world.

There is a left-invariant nearly K\"ahler structure on $S^3 \times S^3$ which arises from a classical construction of 3-symmetric spaces by Ledger and Obata \cite{LO}. It is shown in \cite{Bu} (see also \cite{Bu2}) that this nearly K\"ahler structure is the only one on $S^3 \times S^3$ up to homothety. 
In fact, the proof of this uniqueness result has been the most difficult step in  the classification of homogeneous nearly K\"ahler structures in dimension six. The main tool is the well-known characterisation \cite{RC} of a nearly K\"ahler structure on a six-manifold as an $\SU(3)$-structure $(\omega, \psi^+,\psi^-)$ satisfying the exterior system
\begin{eqnarray}
d\o &=& 3 \psi^+, \label{eq1}\\
d\psi^- &=& \nu \, \o \wedge \o \label{eq2}
\end{eqnarray}
for a real constant $\nu$ which depends on sign and normalisation conventions.

The starting point of this article is the following observation.  The construction of a 3-symmetric space from $G=\SL(2,\bR)$ instead of $\SU(2)$ defines a left-invariant nearly pseudo-K\"ahler structure on $\SL(2,\bR) \times \SL(2,\bR)$. We shortly recall this construction explicitly. The group
$G\times G \times G$ admits a symmetry of order three given by $(g_1,g_2,g_3)
\mapsto (g_2,g_3,g_1)$ which
stabilises the diagonal $\Delta.$ The tangent space of $M^6=G \times G \times G/\Delta$ is
identified with
$$\mathfrak{p}= \{ (X,Y,Z) \in \mathfrak{g}\oplus \mathfrak{g}\oplus\mathfrak{g}\,|\, X+Y+Z=0 \}.$$
Denote by $K_\mathfrak{g}$ the Killing form of $\mathfrak{g}$ and define an
invariant scalar product on $\mathfrak{g}\oplus
\mathfrak{g}\oplus \mathfrak{g}$ by $g=K_\mathfrak{g} \oplus
K_\mathfrak{g}\oplus K_\mathfrak{g}.$ 
This yields a naturally reductive metric on $M^6.$ Using Proposition 5.6 of \cite{G2} this metric is nearly pseudo-K\"ahler.
For completeness sake we recall that the complex structure is given by 
$$ J(X,Y,Z)= \frac{2}{\sqrt 3}(Z,X,Y)+\frac{1}{\sqrt 3}(X,Y,Z).$$

Considering Butruille's results, it is natural to ask how many left-invariant nearly pseudo-K\"ahler structures there are on $\SL(2,\bR) \times \SL(2,\bR)$. Comparing with the results mentioned in the last paragraph, the answer seems a priori hard to guess. The main result of this article is the proof that there is a unique left-invariant nearly pseudo-K\"ahler structure on all Lie groups with Lie algebra $\mathfrak{sl}(2,\bR) \oplus \mathfrak{sl}(2,\bR)$. A byproduct of the proof is the result that there are no nearly para-K\"ahler structures on these Lie groups. We add the remark that there exist co-compact lattices for these Lie groups. Indeed, the article \cite{RV} contains a complete list of the compact quotients of Lie groups with Lie algebras $\mathfrak{sl}(2,\bR)$, which also give rise to compact quotients on a direct product of such groups.

When dealing with nearly pseudo-K\"ahler structures, the problem arises that many facts which are well-known for the Riemannian signature have never been shown for indefinite metrics. For instance, a hyperbolically nearly K\"ahler structure is defined in \cite{B} as a $\SU(2,1)$-structure satisfying the same exterior system \eqref{eq1}, \eqref{eq2} as a nearly K\"ahler structure. It is not obvious whether this definition is equivalent to Gray's classical definition of an (indefinite) nearly pseudo-K\"ahler manifold which is used in \cite{CS,S1}. However, the proof of our main result in section three essentially relies on this exterior system and we have to prove the equivalence.

The close analogy between the pseudo-Hermitian and the para-Hermitian case makes it desirable to give a unified proof dealing with all possible cases at the same time. Therefore, we seize the opportunity and introduce a language that allows us to treat analogous aspects of almost pseudo-Hermitian and almost para-Hermitian geometry simultaneously. This language is consistent with \cite{S3} and similar to \cite{Ki}. In the preliminary section, the necessary basic notions 
are recalled in this unified language. In particular, we recall some facts about stable forms in dimension six which turn out to be very useful in characterising special almost Hermitian structures and special almost para-Hermitian structures. 

Section \ref{NK_section} is devoted to proving the mentioned characterisation of six-dimensional nearly pseudo-K\"ahler and nearly para-K\"ahler manifolds by the exterior system. Since we have to generalise many facts from the Riemannian setting, we give a self-contained proof. Although we follow the ideas of the proof in \cite{RC}, we clarify the structure of the proof by elaborating the role of the Nijenhuis tensor. In particular, we prove that a half-flat structure is additionally nearly half-flat if and only if the Nijenhuis tensor is skew-symmetric. As a first application, we prove some results on the automorphism group of a nearly (para-)K\"ahler six-manifold in section \ref{autom_subsection}.

In section \ref{mainsection}, we finally obtain the aforementioned structure results on $\SL(2,\bR) \times \SL(2,\bR)$. It turns out that the proof is considerably more technical than in the compact case $S^3 \times S^3$, cf. \cite{Bu} or \cite{Bu2}. We also extend the results on $S^3 \times S^3$ by proving the non-existence of nearly (para-)K\"ahler structures of indefinite signature.

The authors wish to thank V. Cort\'es and P.-A. Nagy for useful discussions. In particular, section \ref{autom_subsection} has been inspired by P.-A. Nagy.

\newpage
\section{Preliminaries}

\subsection{Almost pseudo-Hermitian and almost para-Hermitian geometry} 
\label{sect_alm}
We recall that an {\cmssl almost para-complex structure} on a $2m$-dimensional manifold $M$ is an endomorphism field squaring to the identity such that both
eigendistributions (for the eigenvalues $\pm 1$) are $m$-dimensional.
An {\cmssl almost para-Hermitian structure} consists of a neutral metric and an antiorthogonal almost para-complex structure. For a survey on para-complex geometry we refer to \cite{AMT} or \cite{CFG}. 

In the following, we introduce the unified language describing almost pseudo-Hermitian and almost para-Hermitian geometry simultaneously. The philosophy is to put an ``$\ve$'' in front of all notions which is to be replaced by ``para'' for $\ve=1$ and is to be replaced by ``pseudo'' or to be omitted for $\ve = -1$. From now on, we always suppose $\ve \in \{ \pm 1 \}$.

To begin with, we consider the  {\cmssl $\ve$-complex numbers} $\eC = \{ x + \hati y \, , \, x,y \in \mathbb R \}$ with $\hati^2 = \ve$.   
 For the para-complex numbers, $\ve = 1$, there are obvious analogues of conjugation, real and imaginary parts and the square of the  (not necessarily positive) absolute value given by $|z|^2= z  \bar z.$ 

Moreover, let $V$ be a real vector space of even dimension $n=2m$. We call an endomorphism $\pJ$ an {\cmssl $\ve$-complex structure} if $\pJ^2 = \ve {\rm id}_V$ and if additionally, for $\ve=1$,  the $\pm 1$-eigenspaces $V^{\pm}$ are $m$-dimensional. An {\cmssl $\ve$-Hermitian structure} is an $\ve$-complex structure $\pJ$ together with a pseudo-Euclidean scalar-product $g$ which is
{\cmssl $\ve$-Hermitian} in the sense that it holds 
\begin{equation*}
g(\pJ\cdot,\pJ\cdot)=-\ve g(\cdot,\cdot).
\end{equation*}
We denote the stabiliser in $\GL(V)$ of an
$\ve$-Hermitian structure as the {\cmssl $\ve$-unitary group}
$$ \U^\ve(p,q) = \{ L \in \GL(V) \, | \, [L,\pJ] = 0 ,\, L^*g = g \}  \cong
\begin{cases} \U(p,q)  , \, p+q = m , & \mbox{for $\ve = -1$,} \\
\GL(m,\mathbb R) , & \mbox{for $\ve = 1$.}
\end{cases}
$$
Here, the pair $(2p,2q)$ is the signature\footnote{Please note that in our
  convention $2p$ refers to the negative directions.} of the metric for $\ve=-1$. For $\ve=1$, the group $\GL(m,\bR)$ acts reducibly such that $V=V^+ \oplus V^-$ and the signature is always $(m,m)$. 

An {\cmssl almost $\ve$-Hermitian manifold} is a manifold $M$ of dimension $n=2m$ endowed with a $\U^\ve(p,q)$-structure or, equivalently, with an {\cmssl almost $\ve$-Hermitian structure} which consists of an almost $\ve$-Hermitian structure $\pJ$ and an $\ve$-Hermitian metric $g$. The non-degenerate two-form $\omega := g(\cdot , \pJ \cdot )$ is called {\cmssl fundamental two-form}.

Given an almost $\ve$-Hermitian structure $(g,\pJ,\o)$, there exist pseudo-orthonormal local frames $\{ e_1, \dots , e_{2m}\}$ such that $J e_i = e_{i+m}$ for $i=1, \dots, m$ and $ \o = \ve \sum_{i=1}^m \sigma_i e^{i(i+m)}$, 
where $\sigma_i:=g(e_i,e_i)$ for $i=1, \dots, m$. Upper indices will always denote dual (not metric dual) one-forms and $e^{ij}$ stands for $e^i \wedge e^j$. We call such a frame {\cmssl $\ve$-unitary}. If $m \ge 3$, we can always achieve $\sigma_1=\sigma_2$ by reordering the basis vectors.
 
For both values of $\ve$, the $\ve$-complexification $TM \otimes \eC$ of the tangent bundle decomposes into the $\pm \hati$-eigenbundles $TM^{1,0}$ and $TM^{0,1}$. This induces the well-known bi-grading of $\eC$-valued exterior forms 
$$\Omega^{r,s}  = \Gamma (\Lambda^{r,s})= \Gamma (\Lambda^r (TM^{1,0})^* \otimes \Lambda^s (TM^{0,1})^*).$$ 
If $X$ is a vector field on $M$, we use the notation 
$$ X^{1,0} = \frac{1}{2} (X + \hati \ve \pJ X)  \in \Gamma(TM^{1,0}) \: , \:  X^{0,1} = \frac{1}{2} ( X - \hati \ve \pJ X ) \in \Gamma(TM^{0,1}), $$
for the real isomorphisms from $TM$ to $TM^{1,0}$ respectively $TM^{0,1}$. As
usual in almost Hermitian geometry, we define the bundles $\llbracket \Lambda^{r,s} \rrbracket$ for $r\ne s$ and $[\Lambda^{r,r}]$ by the property  
\begin{eqnarray*}
\llbracket \Lambda^{r,s} \rrbracket \otimes \eC =& 
\llbracket \Lambda^{r,s} \rrbracket \oplus \hati \llbracket \Lambda^{r,s} \rrbracket &=
\Lambda^{r,s} \op \Lambda^{s,r} ,\\
{[ \Lambda^{r,r} ]}  \otimes \eC 
=& [ \Lambda^{r,r} ] \oplus \hati [ \Lambda^{r,r} ] &= \Lambda^{r,r}.   
\end{eqnarray*}
The sections in these bundles are denoted as {\cmssl real forms of type $(r,s) + (s,r)$} respectively {\cmssl of type $(r,r)$} and the spaces of sections by $\llbracket \Omega^{r,s} \rrbracket$ respectively by $[\Omega^{r,r}]$ . 
For instance, it holds 
\bea 
[ \Omega^{1,1} ]  &=& \{ \alpha \in \Omega^2 M \, | \, \alpha(X, Y) = \ve \alpha(\pJ X, \pJ Y)\}, \nonumber 
\eea
such that the fundamental form is of type $(1,1)$ and similarly
\bea
\llbracket \Omega^{3,0} \rrbracket  &=& \{ \alpha \in \Omega^3 M \, | \, \alpha(X, Y, Z) = \ve \alpha(X, \pJ Y, \pJ Z)\}. \label{3,0}
\eea

Only in the para-complex case, $\ve = 1$, there is a decomposition of the real tangent bundle 
$TM = \mathcal V \oplus \mathcal H$
into the $\pm 1$-eigenbundles of $\pJ$ which also induces a bi-grading of real forms. It is also straightforward to show that 
\bea 
\llbracket \Lambda^{3,0} \rrbracket \cong \Lambda^3 \mathcal V^* \oplus \Lambda^3 \mathcal H^*, 
\label{para3,0} 
\eea
when considering the characterisation (\ref{3,0}).

Returning to analogies, we recall that the Nijenhuis tensor of the almost $\ve$-complex structure $\pJ$ satisfies
\bea N(X,Y) &=& -\ve [X,Y]-[\pJ X,\pJ Y]+\pJ[\pJ X,Y]+\pJ[X,\pJ Y] \nonumber \\
&=& - \nJ{\pJ X}{Y} + \nJ{\pJ Y}{X} + \pJ \nJ{X}{Y} - \pJ \nJ{Y}{X} \label{N2} \eea
for real vector fields $X,Y,Z$ and for any torsion-free connection $\n$ on
$M$. For both values of $\ve$, it is well-known that the Nijenhuis tensor is the obstruction to the integrability of the almost $\ve$-complex structure.

In the following, let $\n$ always denote the Levi-Civita connection of the metric $g$ of an almost $\ve$-Hermitian manifold. Differentiating the almost $\ve$-complex structure, its square and the fundamental two-form yields for both values of $\ve$ the formulas
\begin{eqnarray}
\nJ{X}{Y} &=& \n_X (\pJ Y) - \pJ (\n_XY), \nonumber  \\
\nJ{X}{\pJ Y} &=& - \pJ \nJ{X}{Y}, \nonumber  \\
\metric{(\n_X\pJ)Y}{Z} &=& -(\n_X\o)(Y,Z),   \label{A0}
\end{eqnarray}
for all vector fields $X,Y,Z$. Using these formulas, it is easy to show that for any \aeh manifold, the tensor $A$ defined by
\begin{equation*}
A(X,Y,Z) = \metric{(\n_X\pJ)Y}{Z} = -(\n_X\o)(Y,Z)
\end{equation*}
has the symmetries 
\bea 
A(X,Y,Z) &=& - A(X,Z,Y), \label{A1}\\
A(X,Y,Z) &=& \ve A(X,\pJ Y, \pJ Z) \label{A2} 
\eea 
for all vector fields $X,Y,Z$. 

The decomposition of the $\U^\ve(p,q)$-representation space of tensors with the same symmetries as $A$ into irreducible components leads to a
classification of \aeh manifolds which is classical for $\U(m)$ \cite{GH}. The
para-complex case for the group $\GL(m,\mathbb R)$ is completely worked out in
\cite{GM}. In \cite{Ki}, the Gray-Hervella classes are generalised to
almost $\ve$-Hermitian structures, which are denoted by generalised almost
Hermitian or $\mathcal{GAH}$ structures there. Analogues of all sixteen
Gray-Hervella classes  are established. These are invariant under the respective group action, but obviously not irreducible for the para-Hermitian case when compared to the decomposition in \cite{GM}. 

Finally, we mention the useful formula
\begin{equation}
\label{tolleFormel} 
2(\n_X\o)(Y,Z) = d\o(X,Y,Z) + \ve d\o(X,\pJ Y, \pJ Z) + \ve g(N(Y,Z),\pJ X)
\end{equation}
holding true for all vector fields $X,Y,Z$ on any almost $\ve$-Hermitian manifold. A short direct proof for $\ve=-1$, $g$ Riemannian, is given in \cite{N}, which also holds literally for pseudo-Riemannian metrics and with sign modifications for $\ve=1$. Alternatively, we refer to \cite{KK} for $\ve =1$. 

\subsection{Stable three-forms in dimension six and $\ve$-complex structures}

We review a construction given in \cite{H2} which associates to a stable three-form $\rho$ on a six-dimensional oriented real vector space $V$ an $\ve$-complex structure on the same vector space. Therefore we recall that
a $k$-form $\rho \in \Lambda^k V^*$ is said to be {\cmssl stable} if its orbit $U$ under $\GL(V)$ is open. Denote by $\kappa$ 
the canonical isomorphism \begin{eqnarray*}  \Lambda^k V^* \cong  \Lambda^{6-k} V \otimes \Lambda^6 V^*.\end{eqnarray*} 
For any three-form $\rho$, one considers $K_\rho \, : \, V \rightarrow V \otimes \Lambda^6 V^*$ defined by 
\begin{eqnarray*}
K_\rho(v):=\kappa((v \lrcorner\, \rho) \wedge \rho)
\end{eqnarray*}
and the quartic invariant 
\begin{eqnarray*}
\lambda (\rho) := \frac{1}{6}{\rm tr}\, ( K_\rho^2 ) \: \in (\Lambda^6 V^*)^{\otimes 2}.
\end{eqnarray*} 
This invariant is different from zero if and only if $\rho$ is stable. Since  $L=\Lambda^6 V^*$ is of dimension one, there exists a 
well-defined notion of positivity and norm in $L\ot L.$  Therefore 
we can, by means of the orientation, associate a volume form $\phi(\rho)$ to a stable three-form $\rho$ by
$$\phi(\rho):= \sqrt{|\lambda(\rho)|}.$$ 
Using this volume we define an endomorphism
\begin{eqnarray*}
J_\rho(v):=\frac{1}{\phi(\rho)}K_\rho(v),
\end{eqnarray*}
which can be proven (cf. \cite{H2},\cite{CLSS}) to be an $\ve$-complex structure, where $\ve$ is the sign of $\lambda(\rho).$ 
For both values of $\ve$, a stable three-form is of type $(3,0) + (0,3)$ with respect to its induced $\ve$-complex structure $J_\rho$ or, in other words, 
$$\Psi_\rho = \rho + \hati J_\rho^* \rho$$ 
is a $(3,0)$-form (where $J_\rho^* \rho (X,Y,Z) = \rho(J_\rho X, J_\rho Y, J_\rho Z )$). Moreover, a stable three-form $\rho$ is non-degenerate in the sense that for $v \in V$ 
\begin{equation}
\label{rhonondeg}
v \hook \rho = 0 \: \Rightarrow \: v = 0
\end{equation}
and the induced volume form satisfies the formula
\begin{equation}
\label{phirho}
\phi (\rho) = \frac{1}{2} \, \pJ_\rho^* \rho \wedge \rho.
\end{equation}

Almost all assertions are straightforward to verify when choosing a basis such that the stable three-form is in the normal form 
\begin{equation}
  \label{normalrho}
\rho = e^{123} + \ve (e^{156} + e^{426} + e^{453})
\end{equation}
which satisfies $\lambda(\rho) = 4 \ve (e^{1 \ldots 6})^{\otimes 2}$, $J_\rho^2 = \ve id_V $ and $J(e_i) = \pm e_{i+3}$ for $i=1,2,3$ where the sign $\pm$ depends on the orientation. 

It is worth mentioning that for every stable three-form in the orbit with $\ve=1$, there is also a basis such that 
\begin{equation}
\label{paranormal}
\rho = e^{123} + e^{456}   
\end{equation}
where  $\{ e_1,e_2,e_3 \}$ and $\{ e_4,e_5,e_6 \}$ span the $\pm 1$-eigenspaces $V^\pm$ of $J_\rho$.

\subsection{Structure reduction of almost $\ve$-Hermitian six-manifolds} \label{struct_red_ect}

Let $(V,g,\pJ,\omega)$ be a $2m$-dimensional $\ve$-Hermitian vector space and $\Psi = \psi^+ + \hati \psi^-$ be an $(m,0)$-form of \emph{non-zero} length. We define the special $\ve$-unitary group $\SU^\ve(p,q)$ as the stabiliser of $\Psi$ in the  $\ve$-unitary group $\U^\ve(p,q)$ such that 
$$ \SU^\ve(p,q) = \Stab_{\GL(V)}(g,\pJ,\Psi) \cong
\begin{cases} \SU(p,q) \: , \, p+q = m, & \mbox{for $\ve = -1$}, \\
\SL(m,\mathbb R), & \mbox{for $\ve = 1$,}
\end{cases}
$$
where $\SL(m,\bR)$ acts reducibly such that $V=V^+ \oplus V^-$.

With this notation, an {\cmssl $\SU^\ve(p,q)$-structure} on a manifold $M^{2m}$ is an almost $\ve$-Hermitian structure $(g,\pJ,\o)$ together with a global $(m,0)$-form $\Psi$ of non-zero constant length. Locally, there exists an $\ve$-unitary frame $\{ e_1, \dots,e_m,e_{m+1}=Je_1,\dots , e_{2m}=Je_m\}$ which is adapted to the $\SU^\ve(p,q)$-reduction $\Psi$ in the sense that 
\begin{equation}
\label{adapted}
\Psi = a (e^1 + \hati e^{(m+1)}) \wedge .\,.\,. \wedge (e^m + \hati e^{2m})
\end{equation}
for a constant $a \in \bR^*$. 

In dimension six, there is a characterisation of $\SU^\ve(p,q)$-structures in terms of stable forms. Given a six-dimensional real vector space $V,$ we call a pair $(\omega,\rho)$ of a stable  $\omega \in \Lambda^2 V^*$ and a stable  $\rho \in \Lambda^3 V^*$ compatible if it holds
\begin{eqnarray}
 \label{comp1} \omega \wedge \rho &=&0.
\end{eqnarray}
We claim that the stabiliser in $\GL(V)$ of a compatible pair is 
$$ \Stab_{\GL(V)}(\omega,\rho) = \SU^\ve(p,q) , \quad p + q =3, $$
where $\varepsilon \in \{ \pm 1 \}$ is the sign of $\lambda(\rho)$, that is
$J_\rho^2 = \varepsilon id_V$. This can be seen as follows. For the two-form $\omega$, stability is equivalent to non-degeneracy and we choose the orientation on $V$ such that $\omega^3$ is positive. By the previous section, we can associate an $\ve$-complex structure $J_\rho$ to the stable three-form $\rho$. For instance in an adequate basis, it is easy to verify that $\omega \wedge \rho = 0$ is equivalent to the skew-symmetry of $J_\rho$ with respect to $\omega$. Equivalently, the pseudo-Euclidean metric 
\begin{equation} \label{inducedmetric}
g = \varepsilon \, \omega(\cdot,J_\rho\cdot),
\end{equation}
induced by $\omega$ and $\rho$ is $\ve$-Hermitian with respect to $J_\rho$. Since $\Psi_\rho = \rho + \hati J_\rho^* \rho$ is a $(3,0)$-form and the stabiliser of $\omega$ and $\rho$ also stabilises the tensors induced by them, the claim follows.

We conclude that an $\SU^\ve(p,q)$-structure, $p+q =3$, on a six-manifold is
characterised by a pair $(\o,\psi^+) \in \Omega^2 M \times \Omega^3 M$ of everywhere stable and compatible forms such that the induced $(3,0)$-form $\Psi = \psi^+ + \hati J_{\psi^+}^* \psi^+ = \psi^+ + \hati \psi^-$ has constant non-zero length with respect to the induced metric \eqref{inducedmetric}. 
In an $\ve$-unitary frame which is adapted to $\Psi$ in the sense of \eqref{adapted}, the formula
\begin{equation}
\label{lengthpsi}
\psi^- \wedge \psi^+ = \| \psi^+ \|^2 \frac{1}{6}\, \omega^3  
\end{equation}
is easily verified. Thus, given a compatible pair $(\o,\psi^+)$ of stable forms, it can be checked that the induced $(3,0)$-form $\Psi$ has constant non-zero length without explicitly computing the induced metric.

We remark that in the almost Hermitian case, the literature often requires $\Psi$ to be normalised such that $\| \psi^+ \|^2 = 4$, for instance in \cite{ChSa}. In the more general almost $\ve$-Hermitian case, we have in an adapted local $\ve$-unitary frame \eqref{adapted} with $\sigma_1 = \sigma_2$
\begin{equation}
\label{adaptedpsi}
\psi^+ = a (e^{123} + \ve (e^{156} + e^{426} + e^{453})) \qquad \mbox{and} \qquad  \| \psi^+ \|^2 = 4 a^2 \sigma_3
\end{equation}
for a real constant $a$. Therefore we have to consider two different
normalisations  $\| \psi^+ \| = \pm 4$ or we multiply the metric by $-1$ if necessary such that $\| \psi^+ \|$ is always positive.

Finally, we remark that $\SU(3)$-structures are classified in \cite{ChSa} and it is shown that the intrinsic torsion is completely determined by the exterior derivatives $d\omega$, $d\psi^+$ and $d\psi^-$.

\section{Nearly pseudo-K\"ahler and nearly para-K\"ahler manifolds}
\label{NK_section}

The main objective of this section is to generalise the characterisation of
six-dimensional nearly K\"ahler manifolds by an exterior differential system to nearly pseudo-K\"ahler and nearly para-K\"ahler manifolds. 

\subsection{General properties}
\bd \label{NK_Def}
An \aeh manifold $(M^{2m},g,\pJ,\o)$ is called {\cmssl nearly $\ve$-K\"ahler manifold}, if its Levi-Civita  connection $\n$ satisfies the {\cmssl nearly $\ve$-K\"ahler condition} 
\be (\n_X\pJ)\,X = 0, \quad \forall   X \in \Gamma(TM). \nonumber \ee
A nearly $\ve$-K\"ahler manifold is called {\cmssl strict} if $\n_X \pJ \ne 0$ for all non-trivial vector fields $X$.
\ed
A tensor field $B \in \Gamma ((TM^*)^{\otimes 2} \otimes \, TM)$ is called totally skew-symmetric if the tensor $g(B(X,Y),Z)$ is a three-form. The following characterisation of a nearly $\ve$-K\"ahler manifold is well-known in the Riemannian context.
\bp 
\label{char_NK}
An \aeh manifold $(M^{2m},g,\pJ,\o)$ satisfies the nearly $\ve$-K\"ahler condition if and only if 
$d\o$ is of real type $(3,0) + (0,3)$ and the Nijenhuis tensor is totally skew-symmetric. 
\ep
\begin{proof}
The nearly $\ve$-K\"ahler condition is satisfied if and only if the tensor $A = - \nabla \omega$ is a three-form because of the antisymmetry (\ref{A1}). 

Assume first that $(g,\pJ,\o)$ is a nearly $\ve$-K\"ahler structure. Comparing the identities (\ref{3,0}) and (\ref{A2}), we see that the real three-form $A$ is of type $(3,0) + (0,3)$. Since $d\o$ is the alternation of $\n\o$, we have 
\begin{equation}
\label{dovsA}
d\o =3\n \o=-3A \in \llbracket \Omega^{3,0} \rrbracket.  
\end{equation}
Furthermore, if we apply the nearly $\ve$-K\"ahler condition to the expression (\ref{N2}), the Nijenhuis tensor of a nearly $\ve$-K\"ahler structure simplifies to 
\begin{equation}
\label{NK_N}
N(X,Y) = 4 \, \pJ \nJ{X}{Y}.
\end{equation}
We conclude that the Nijenhuis tensor is skew-symmetric since  
\begin{equation}
\label{NvsA}
 g(N(X,Y),Z) = - 4 A(X,Y,\pJ Z) \stackrel{\eqref{A2}}{=} - 4 \ve {\pJ}^*A( X, Y, Z).
\end{equation}

The converse follows immediately from the identity \eqref{tolleFormel} when considering \eqref{3,0}. For self-containedness we give a direct proof. Assume that $d\o \in \llbracket \Omega^{3,0} \rrbracket $ and the Nijenhuis tensor is skew-symmetric. To begin with, we observe that 
$$ \no{Y}{X}{X} = 0 = \no{JY}{X}{\pJ X}$$
by \eqref{A1} and \eqref{A2}. With this identity, we have on the one hand 
\begin{eqnarray*}
 0 &=& \ve g(N(\pJ X,\pJ Y),\pJ X) = g(N(X,Y),\pJ X) \\ 
&\stackrel{\eqref{N2}}{=}& -g(\nJ{\pJ X}{Y},\pJ X) + g(\nJ{\pJ Y}{X},\pJ X) + g(\pJ \nJ{X}{Y},\pJ X) - g(\pJ \nJ{Y}{X}, \pJ X) \\
&\stackrel{\eqref{A0}}{=}& \no{\pJ X}{Y}{\pJ X} + \ve \no{X}{Y}{X} \\
&\stackrel{\eqref{A1}}{=}& \no{\pJ X}{Y}{\pJ X} - \ve \no{X}{X}{Y} ,
\end{eqnarray*}
and on the other hand
\begin{eqnarray*}
0 &=& \ve d\o(X,X,Y) \stackrel{\eqref{3,0}}{=} d\o(X, \pJ X,\pJ Y) \\
&=& (\n_{X} \o)({\pJ X},\pJ Y) + (\n_{\pJ X} \o)(\pJ Y,X) + (\n_{\pJ Y} \o)(X,{\pJ X})\\
&\stackrel{\eqref{A2}}{=}& \ve (\n_X \o)({X},{Y}) + (\n_{\pJ X} \o)(Y,{\pJ X}).
\end{eqnarray*} 
It follows that $\no{X}{X}{Y}=0$ which is equivalent to the nearly $\ve$-K\"ahler condition.
\end{proof}
\br
The notion of nearly $\ve$-K\"ahler manifold corresponds to the generalised class $\mathcal W_1$ in \cite{Ki}. However, in the para-Hermitian case, there are two subclasses, see \cite{GM}. Indeed, we already observed that 
$$A = -\n \omega \in \llbracket \Omega^{3,0} \rrbracket \stackrel{\eqref{para3,0}}{=} \Gamma (\Lambda^3 \mathcal V^* \oplus \Lambda^3 \mathcal H^*) $$
for a nearly para-K\"ahler manifold. 
\er

We call a connection $\bn$ on an almost $\ve$-Hermitian manifold $(M^{2m},g,\pJ,\o)$ {\cmssl $\ve$-Hermitian} if $\bn g = 0$ and $\bn J =0$. 
\begin{Prop}
An almost $\ve$-Hermitian manifold $(M^{2m},g,\pJ,\o)$ admits an $\ve$-Hermitian connection with totally skew-symmetric torsion if and only if the Nijenhuis tensor is totally skew-symmetric. If this is the case, the connection $\bn$ and its torsion $T$ are uniquely defined by 
\begin{eqnarray}
g(\bn_X Y,Z) &=&  g(\n_X Y,Z) + \frac{1}{2} g(T(X,Y),Z), \nonumber \\
g(T(X,Y),Z) &=& \ve g(N(X,Y),Z) -d\o(JX,JY,JZ), \nonumber
\end{eqnarray}
and we call $\bn$ the canonical $\ve$-Hermitian connection (with skew-symmetric torsion).
\end{Prop}
\begin{proof}
The Riemannian case is proved in \cite{FI}, the para-complex case in \cite{IZ}. In fact, the sketched proof in \cite{FI} holds literally for the almost pseudo-Hermitian case with indefinite signature as well. For completeness, we give a direct proof for all cases simultaneously.

Let $T(X,Y) = \bn_X Y - \bn_Y X - [X,Y] = S_X Y - S_Y X$ be the totally skew-symmetric torsion of an $\ve$-Hermitian connection $\bn$ where $S_XY=\bar \n_XY -\n_XY$ is the difference tensor with respect to the Levi-Civita connection $\n$ of $g$. Then, the Nijenhuis tensor is totally skew-symmetric as well, since we have
\begin{equation}
\label{NT}
g(N(X,Y),Z) = \ve g(T(X,Y),Z) + g(T(JX,JY),Z) +g(T(JX,Y),JZ) +g(T(X,JY),JZ),
\end{equation}
using only $\bn J=0$. Moreover, the difference tensor $S_X$ is skew-symmetric with respect to $g$, for $\bn g=0$. Combining this fact with the total skew-symmetry of the torsion, cf. for example \cite{CS} Lemma 1, 
we find that $S_XY = - S_Y X$ and consequently 
\begin{equation*}
g(\bn_X Y,Z) =  g(\n_X Y,Z) + \frac{1}{2} g(T(X,Y),Z).  
\end{equation*}
With this identity and $\bn \o=0$, the equation
\begin{equation}
\label{no}
2 \n_{JX}\o (Y,Z) 
= g(T(JX,Y),JZ) + g(T(JX,JY),Z)
\end{equation}
follows. Finally, we verify the claimed formula for the torsion:
\begin{eqnarray*}
d\o(JX,JY,JZ) 
&\stackrel{\eqref{A2}}{=}&
\ve ( \n_{JX}\o (Y,Z) +  \n_{JY}\o (Z,X)  + \n_{JZ}\o (X,Y) ) \\
&\stackrel{\eqref{no}}{=}&
\ve( g(T(JX,JY),Z) + g(T(JX,Y),JZ) + g(T(X,JY),JZ) ) \\
&\stackrel{\eqref{NT}}{=}& 
\ve g(N(X,Y),Z) - g(T(X,Y),Z).
\end{eqnarray*}
Conversely, if the Nijenhuis tensor is skew-symmetric, is is straightforward to verify that the defined connection is $\ve$-Hermitian with skew-symmetric torsion.
\end{proof}
\begin{remark}
An almost Hermitian manifold is said to be of type $\mathcal G_1$ if it admits
a Hermitian connection with skew-symmetric torsion, see for example \cite{N}.
More generally, the proposition justifies to say that an almost $\ve$-Hermitian manifold is of type $\mathcal G_1$ if it admits an $\ve$-Hermitian connection with skew-symmetric torsion.
\end{remark}
In particular, the proposition applies to nearly $\ve$-K\"ahler manifolds. In this case, the skew-symmetric torsion $T$ of the canonical $\ve$-Hermitian connection simplifies to 
$$ T(X,Y) = \varepsilon \pJ(\n_X \pJ)Y = \frac{1}{4} \varepsilon N(X,Y)$$
due to the identities \eqref{dovsA}, \eqref{NK_N} and \eqref{NvsA}.
\bp 
\label{paralleltorsion}
The canonical $\ve$-Hermitian connection $\bar \n$ of a nearly $\ve$-K\"ahler manifold $(M^{2m},\pJ,g,\omega)$ satisfies $$\bar\n (\n \pJ)=0 \quad \mbox{and} \quad \bar\n(T)=0.$$
\ep
\begin{proof} 
The two assertions are equivalent since $\bar \n \pJ =0$. A short proof of the first assertion for the Hermitian case is given in \cite{BM}. This proof generalises without changes to the pseudo-Hermitian case since it essentially uses the identity
\bean \label{id_equ} 2g((\n^2_{W,X}J)Y,Z)= -\sigma_{X,Y,Z} \, g((\n_WJ)X,(\n_YJ)JZ), \eean
which was proved in \cite{G1} for Riemannian metrics and also holds true in the pseudo-Riemannian setting (\cite[Proposition 7.1]{Ka}). The para-Hermitian version is proved in \cite[Theorem 5.3]{IZ}.
\end{proof}
\begin{Cor}
\label{Aconst}
On a nearly $\ve$-K\"ahler manifold $(M^{2m},\pJ,g,\omega)$, the tensors $\n J$ and $N=4\ve T$ have constant length.
\end{Cor}
\begin{proof}
This is obvious since both tensors are parallel with respect to the connection $\bar \n$ which preserves in particular the metric.
\end{proof}
\br 
In dimension six, the fact that $\n J$ has constant length is usually expressed by the equivalent assertion that a nearly $\ve$-K\"ahler six-manifold is of constant type, i.\ e.\ there is a constant $\kappa \in \bR$ such that 
\bean
g(\nJ{X}{Y},\nJ{X}{Y}) = \kappa\,\{\, g(X,X) g(Y,Y) - g (X,Y)^2 + \ve g(\pJ X, Y)^2 \,\}. 
\eean
In fact, the constant is $\kappa = \frac{1}{4} \| \n J \|^2 $. 
Furthermore, it is well-known in the Riemannian case that strict nearly K\"ahler six-manifolds are Einstein manifolds with Einstein constant $5 \kappa$ \cite{G1}. The same is true in the para-Hermitian case \cite{IZ} and in the pseudo-Hermitian case \cite{S4}.
\er
The case $\| \n J \|^2 =0 $ for a strict nearly $\ve$-K\"ahler six-manifold can only occur in the para-complex world. We give different characterisations of such structures which provide an obvious break in the analogy of nearly para-K\"ahler and nearly pseudo-K\"ahler manifolds.

\bp \label{const_type_zero_prop} 
For a six-dimensional strict nearly para-K\"ahler manifold $(M^6,g,J,\omega)$, the following properties are equivalent:
\begin{enumerate}
\item [(i)]  $\| \n \pJ \|^2 = \| A \|^2 = 0$ 
\item [(ii)] The three-form $A=-\n \omega \in \llbracket \Omega^{3,0} \rrbracket$ is either in $\Gamma (\Lambda^3 \mathcal V^*)$ or in $\Gamma (\Lambda^3 \mathcal H^*)$.
\item [(iii)] The three-form $A=-\n \omega \in \llbracket \Omega^{3,0} \rrbracket$ is not stable.
\item [(iv)] The metric $g$ is Ricci-flat.
\end{enumerate}
\ep
\begin{proof}
We choose a local frame $\{ e_1, \dots ,e_6 \}$ such that $\{ e^1, e^2, e^3 \}$ spans the $+1$-eigenspace $\mathcal V^*$ of $J$, $\{ e^4, e^5, e^6 \}$ spans the  $-1$-eigenspace $\mathcal H^*$ of $J$ and $g(e_i,e_{i+3})=1$ for $i=1,2,3$. According to \eqref{para3,0}, there are local functions $a$ and $b$ such that $A = a e^{123} + b e^{456}$. Thus, it holds 
$$ \| A \|^2 \frac{1}{6} \omega^3  \stackrel{\eqref{lengthpsi}}{=} \pJ^* A \wedge A = (a e^{123} - b e^{456}) \wedge (a e^{123} + b e^{456}) = 2 \, a \, b \, e^{123456} .$$
With $\omega = - e^{14} - e^{25} - e^{36}$ and $\omega^3 = 6 e^{123456}$, we have $\| A \|^2 = 2 a b$. Since $A$ is nowhere zero due to the strictness and considering also \eqref{phirho}, the first three assertions are equivalent to $a=0$ or $b=0$. Finally, assertions (i) and (iv) are equivalent by \cite[Theorem 5.5]{IZ}.
\end{proof}

Flat strict nearly para-K\"ahler manifolds $(M,g,J,\omega)$ are classified in \cite{CS2}. It turns out that they always satisfy  $\| \n \pJ \|^2 = 0$. 
In \cite{GM}, almost para-Hermitian structures on tangent bundles $TN$ of 
real three-dimensional manifolds $N^3$ are discussed. It is shown that the existence
of nearly para-K\"ahler manifolds satisfying the second condition of Proposition
\ref{const_type_zero_prop} is equivalent to the existence of a certain 
connection on $N^3.$ However, to the authors best knowledge, 
there exists no reference for an example of a Ricci-flat nearly para-K\"ahler structure which is not flat.

\subsection{Characterisations by exterior differential systems in dimension six} 

The following lemma explicitly relates the Nijenhuis tensor to the exterior differential.
For $\ve = -1$, it gives a characterisation of Bryant's notion of a quasi-integrable $\U(p,q)$-structure, $p+q = 3$, in dimension six \cite{B}. 

Let $(M^6,g,\pJ,\o)$ be a six-dimensional almost $\ve$-Hermitian  manifold. If $\{ e_1,\dots, e_6=J e_3 \}$ is a local $\ve$-unitary frame, we define a local frame $\{ E^1, E^2, E^3 \}$ of $(TM^{1,0})^*$ by 
$$E^i := (e^i + \hati \ve \pJ e^i) = (e^i + \hati e^{i+m})$$
for $i = 1,2,3$ and denote it as a local $\ve$-unitary frame of $(1,0)$-forms. The dual vector fields of the $(1,0)$-forms are 
$$ E_i = \, e_i^{1,0} = \frac{1}{2} (e_i +  \hati \ve \pJ e_i) = \frac{1}{2} (e_i + \hati \ve e_{i+m}) $$
such that the $\eC$-bilinearly extended metric satisfies 
$$ g(E_i,\bar E_j) = \frac{1}{2} \sigma_i \delta_{ij} \qquad \mbox{and} \qquad g(E_i,E_j) = 0$$
in such a frame.
\bl 
\label{lambdalemma}
The Nijenhuis tensor of an almost $\ve$-Hermitian six-manifold $(M^6,g,\pJ,\o)$ is totally skew-symmetric if and only if for every local $\ve$-unitary frame of $(1,0)$-forms, there exists a local $\eC$-valued function $\lambda$ such that
\begin{eqnarray} 
\label{Nij_tot_skew_equ} 
(dE^{\tau(1)})^{0,2} =  \lambda\,\sigma_{\tau (1)} \, E^{\overline { \tau(2)}\, \overline{ \tau(3)}}
\end{eqnarray}
for all even permutations $\tau$ of $\{ 1, 2, 3 \}$.
\el
\pf 
First of all, the identities 
\bean N(\bar V,\bar W) &=& -4 \ve [\bar V,\bar W]^{1,0} \quad \mbox{and} \quad N(V, \bar W) = 0 \eean
for any vector fields $V=V^{1,0}$, $W=W^{1,0}$ in $TM^{1,0}$ follow
immediately from the definition of $N$. 
Using the first identity, we compute in an arbitrary local $\ve$-unitary frame 
\bea 
\label{Nij_tot_help} 
d E^i(\bar E_j,\bar E_k) &=& - E^i( [\bar E_j,\bar E_k] ) = - 2 \sigma_i \, g( [\bar E_j,\bar E_k],\bar E_i) \nonumber \\
&=& - 2 \sigma_i \, g( [\bar E_j,\bar E_k]^{1,0},\bar E_i) = \frac{1}{2} \ve
\,\sigma_i \, g (N(\bar E_j,\bar E_k),\bar E_i) \nonumber
\eea
for all possible indices $1\le i,j,k \le 3$. 
If the Nijenhuis tensor is totally skew-symmetric, equation \eqref{Nij_tot_skew_equ} follows by setting 
\begin{equation}
\label{lambda}
\lambda = \frac{1}{2} \ve  \, g (N(\bar E_1,\bar E_2),\bar E_3).
\end{equation}
Conversely, the assumption \eqref{Nij_tot_skew_equ} for every local $\ve$-unitary frame implies that the Nijenhuis tensor is everywhere a three-form when considering the same computation and $N(V, \bar W) = 0$.
\epf

If there is an $\SU^\ve(p,q)$-reduction with closed real part, this characterisation can be reformulated globally in the following sense.

\bp \label{zweiteGleichung} 
Let $(\o,\psi^+)$ be an $\SU^\ve(p,q)$-structure on a six-manifold $M$ such that $\psi^+$ is closed. Then the Nijenhuis tensor is totally skew-symmetric if and only if
\begin{equation}
 d\psi^- = \nu \,\o \wedge \o \label{2ndNK}
\end{equation}
for a global real function $\nu$.
\ep
\pf It suffices to proof this locally. 
Let $\{ E^i \}$ be an $\ve$-unitary frame of $(1,0)$-forms with $\sigma_1=\sigma_2$ which is adapted to the $\SU^\ve(p,q)$-reduction such that $\Psi = \psi^+ + \hati \psi^- = a E^{123}$ for a real constant $a$ as in \eqref{adapted}. The fundamental two-form is 
\bean 
\o= - \frac{1}{2} \hati \sum_{k=1}^m \sigma_k \, E^{k \bar k}
\eean
in such a frame.
Furthermore, as $\psi^+$ is closed, we have $d\Psi = \hati d \psi^- = -d\bar \Psi,$ which
implies that $d\Psi \in \Lambda^{2,2}$. Considering this, we compute the real 4-form
\bean d \psi^-  =  \ve \hati \, d \Psi = \ve \hati a \, 
\left( (dE^1)^{0,2} \wedge E^{23} + (dE^2)^{0,2} \wedge E^{31} + (dE^3)^{0,2} \wedge E^{12} \right)
\eean
and compare this expression with
\bean
\o \wedge \o &=& \frac{1}{2} \ve 
( \sigma_2 \sigma_3 \, E^{2\bar 2 3\bar3} 
+ \sigma_1 \sigma_3 \, E^{1\bar 1 3\bar3}
+ \sigma_1 \sigma_2 \, E^{1\bar 1 2\bar2} ) \\
&=& - \frac{1}{2} \ve \sigma_3 
( \sigma_1\, E^{\bar 2 \bar 3 2 3} 
+ \sigma_2\, E^{\bar 3 \bar 1 3 1} 
+ \sigma_3\, E^{\bar 1 \bar 2 1 2}).
\eean
Hence, by Lemma \ref{lambdalemma}, the Nijenhuis tensor is totally
skew-symmetric if and only if $d\psi^- = \nu \,\o \wedge \o$ holds true for a real function $\nu$. More precisely, the two functions $\nu$ and $\lambda$ are related by the formula 
\begin{equation}
\label{nu}
\nu = - 2 \sigma_3 \hati a \lambda.   
\end{equation}
\epf
An $\SU^\ve(p,q)$-structure $(\o,\psi)$ is called \emph{half-flat} if 
$$ d\psi=0, \quad d \omega^2 = 0,$$
and \emph{nearly half-flat} if 
$$d\psi = \nu \,\o \wedge \o$$ 
for a real constant $\nu$. These notions are defined for the Riemannian signature in \cite{ChSa} respectively \cite{FIMU} and extended to all signatures in \cite{CLSS}.
\bc \label{halfflatplusx}
Let $(\o,\psi^+)$ be a half-flat $\SU^\ve(p,q)$-structure on a six-manifold $M$. Then, the Nijenhuis tensor is totally skew-symmetric if and only if $(\o,\psi^-)$ is nearly half-flat.
\ec
\begin{proof}
If $(\o,\psi^-)$ is nearly half-flat, the equation \eqref{2ndNK} is satisfied
by definition and the Nijenhuis tensor is skew-symmetric by the previous
proposition. In particular one has $d\o^2=0.$ Conversely, if the Nijenhuis tensor is skew, we know that \eqref{2ndNK} holds true for a real function $\nu$, since we have $d\psi^+=0$. Differentiating this equation and using $d\omega^2=0$, we obtain $d \nu \wedge \omega^2 = 0$. The assertion follows as wedging by $\omega^2$ is injective on one-forms.
\end{proof}
\br
An interesting property of $\SU(p,q)^\ve$-structures which are both half-flat and nearly half-flat in the sense of the corollary is the fact that, given that the manifold and the $\SU(p,q)^\ve$-structure are analytic, the structure can be evolved to both a parallel $G_2$-structure and a nearly parallel $G_2$-structure via the Hitchin flow. For details, we refer to \cite{H1} and \cite{St} for the compact Riemannian case and \cite{CLSS} for the non-compact case and indefinite signatures.

In \cite{ChSw}, six-dimensional nilmanifolds $N$ admitting an invariant
half-flat $\SU(3)$-structure $(\o,\psi^+)$ such that $(\o,\psi^-)$ is nearly
half-flat are classified. As six nilmanifolds admit such a structure, we
conclude that these structures are not as scarce as nearly K\"ahler
manifolds. It is also shown in this reference, that these structures induce invariant $G_2$-structures with torsion on $N \times S^1$.

We give another example of a (normalised) left-invariant $\SU(3)$-structure on $S^3 \times S^3$ which satisfies $d\psi^+ = 0, \; d\psi^- = \,\o \wedge \o$ such that $d \omega$ neither vanishes nor is of type (3,0) + (0,3). We choose a global frame of left-invariant vector fields $\{ e_1, e_2, e_3, f_1, f_2, f_3 \}$ on $S^3 \times S^3$ such that
\begin{eqnarray*}
&de^1= e^{23}\, , \quad de^2= e^{31}\, , \quad de^3= e^{12} \quad\, , \quad
df^1= f^{23}\, , \quad df^2= f^{31}\, , \quad df^3= f^{12},&
\end{eqnarray*}
and set with $x= 2 + \sqrt 3$
\begin{eqnarray*}
\omega &=& e^{1}f^{1}+e^{2}f^{2}+e^{3}f^{3}, \\
\psi^+ &=& -\frac{1}{2}x^2 e^{123} +2 x e^{12}f^{3}-2 x e^{13}f^{2}-2 x e^{1}f^{23} 
+ \, 2 x e^{23}f^{1} \\
&& +\, 2 x e^{2}f^{13}-2 x e^{3}f^{12} + (4 x-8) f^{123} ,\\
\psi^- &=& \frac{1}{2} x e^{123}-2  e^{1}f^{23}+2  e^{2}f^{13}-2  e^{3}f^{12}+4  f^{123} ,\\
g &=& x \, ( e^1 )^2 + x \,  ( e^2 )^2 + x \, ( e^3 )^2 + 4 \, ( f^1 )^2 + 4 \, ( f^2 )^2 + 4 \, ( f^3 )^2 \\
&&-\, 2 x \, e^1 \!\cdot\! f^1 - 2 x \, e^2 \!\cdot\! f^2 - 2 x e^3 \!\cdot\! f^3.
\end{eqnarray*}
\er

Finally, we come to the characterisation of six-dimensional nearly $\ve$-K\"ahler manifolds by an exterior differential system generalising the classical result of \cite{RC} which holds for $\ve =-1$ and Riemannian metrics.

\bt 
\label{ext_diff_syst_th} 
Let $(M, g, \pJ, \omega)$ be an almost $\ve$-Hermitian six-manifold. Then $M$ is a strict nearly $\ve$-K\"ahler manifold with $\| \n J \|^2 \ne 0$ if and only if there is a reduction $\Psi = \psi^+ + \hati \psi^-$ to $\SU(p,q)^\ve$ which satisfies
\begin{eqnarray}
d\o &=& 3 \, \psi^+,  \\
d\psi^- &=& 2 \, \kappa  \, \o \wedge \o,
\end{eqnarray}
where $\kappa = \frac{1}{4} \| \n J \|^2 $ is constant and non-zero.
\et
\begin{remark}
Due to our sign convention $\omega=g(.,\pJ.)$, the constant $\kappa$ is positive in the Riemannian case and the second equation differs from that of other authors. Furthermore, we will sometimes use the term nearly $\ve$-K\"ahler manifold of non-zero type if $\| \n J \|^2 \ne 0$.
\end{remark}
\begin{proof}
By Proposition \ref{char_NK}, the manifold 
$M$ is nearly $\ve$-K\"ahler if and only if $d \o$ is of type $(3,0)+(0,3)$ and the Nijenhuis tensor is totally skew-symmetric.

Therefore, when $(g,\pJ,\o)$ is a strict nearly $\ve$-K\"ahler structure such that $ \| A \|^2 = \| \n J \|^2$ is constant (by Corollary \ref{Aconst}) and not zero (by assumption), we can define the reduction $\Psi = \psi^+ + \hati \psi^-$ by $\psi^+ = \frac{1}{3} d\omega = - A$ and $\psi^- = \pJ^* \psi^+$ such that the first equation is satisfied. Since $\omega$ is of type $(1,1)$ and therefore $d (\o \wedge \o) = 2 d\o \wedge \o = 0$, this reduction is half-flat.
Thus, Corollary \ref{halfflatplusx} and the skew-symmetry of $N$ imply that there is a constant $\nu \in \bR$ such that $d\psi^- = \nu \,\o \wedge \o$. 

According to \eqref{adapted}, we can choose an $\ve$-unitary local frame with $\sigma_1=\sigma_2$, such that 
$$ \Psi = - A - \hati \pJ^* A = a E^{123}, $$ 
where $a$ is constant and satisfies $4 \kappa = \| \n \pJ \| ^2 = \| \psi^+ \|^2 = 4 a^2 \sigma_3$ by \eqref{adaptedpsi}. Now, the functions defined in Lemma \ref{lambdalemma} and Proposition \ref{zweiteGleichung} evaluate as
\begin{eqnarray*}
\lambda &\stackrel{\eqref{lambda}}{=}& \frac{1}{2}\ve g( N(\bar E_1,\bar E_2), \bar E_3) \stackrel{\eqref{NvsA}}{=} - 2 J^*A(\bar E_1,\bar E_2,\bar E_3) = - \ve \, \hati a, \\
\nu &\stackrel{\eqref{nu}}{=}& - 2 \sigma_3 \hati a \lambda = 2 \sigma_3 a^2 = 2   \kappa .
\end{eqnarray*}

Conversely, if a given $\SU(p,q)^\ve$-structure satisfies the exterior system,
the real three-form $\psi^+$ is obviously closed and the Nijenhuis tensor is
totally skew-symmetric by Corollary \ref{halfflatplusx}. Considering that
$d \o = 3 \n \o$ is of type $(3,0)+(0,3)$ by the first equation, the structure is nearly $\ve$-K\"ahler. Since $A=-\psi^+$ is stable, the structure is strict nearly $\ve$-K\"ahler by \eqref{rhonondeg} and $\| \n J \| = \| A\| \ne 0$ by Proposition \ref{const_type_zero_prop}. Now, the computation of the constants in the adapted $\ve$-unitary frame shows that in fact $\| \n J \| = 4 \kappa$.
\end{proof}

\subsection{Consequences for automorphism groups} \label{autom_subsection}

An automorphism of an $\SU^\ve(p,q)$-struc-\\ture on a six-manifold $M$ is an automorphism of principal fibre bundles or equivalently, a diffeomorphism of $M$ preserving all tensors defining the $\SU^\ve(p,q)$-structure. By our discussion on stable forms in section \ref{struct_red_ect}, an $\SU^\ve(p,q)$-structure is characterised by a pair of compatible stable forms $(\omega,\rho) \in \Omega^2 M \times \Omega^3 M$. Since the construction of the remaining tensors $J,\psi^-$ and $g$ is invariant, an diffeomorphism preserving the two stable forms is already an automorphism of the $\SU^\ve(p,q)$-structure and in particular an isometry.

This easy observation has the following consequences when combined with the exterior systems of the previous section and the naturality of the exterior derivative.
\bp
Let $(\omega,\psi^+)$ be an $\SU^\ve(p,q)$-structure on a six-manifold $M$.
\begin{enumerate}[(i)]
\item If the exterior differential equation
\[ d\omega = \mu \; \psi^+\] 
is satisfied for a constant $\mu \ne 0$, then a diffeomorphism $\Phi$ of $M$ preserving $\omega$ is an automorphism of the $\SU^\ve(p,q)$-structure and in particular an isometry.
\item If the exterior differential equation 
\[ d\psi^- = \nu \; \omega \wedge \omega \] 
is satisfied for a constant $\nu\ne 0$, then a diffeomorphism $\Phi$ of $M$ preserving
\begin{enumerate}[(a)]
\item the real volume form and $\psi^+,$ 
\item or the real volume form and $\psi^-,$
\item or the $\ve$-complex volume form $\Psi = \psi^+ + \hati \psi^-,$
\end{enumerate}
is an automorphism of the $\SU^\ve(p,q)$-structure and in particular an isometry.
\end{enumerate}
\ep
We like to emphasise that both parts of the Proposition apply to strict nearly $\ve$-K\"ahler structures of non-zero type.

Conversely, it is known for \emph{complete} Riemannian nearly K\"ahler manifolds, that orientation-preserving isometries are automorphism of the almost Hermitian structure except for the round sphere $S^6$, see for instance \cite[Proposition 4.1]{Bu2} in this handbook. However, this is not true if the metric is incomplete. In \cite[Theorem 3.6]{FIMU}, a nearly K\"ahler structure is constructed on the incomplete sine-cone over a Sasaki-Einstein five-manifold $(N^5,\eta,\o_1,\o_2,\o_3)$. In fact, the Reeb vector field dual to the one-form $\eta$ is a Killing vector field which does not preserve $\o_2$ and $\o_3$. Thus, by the formulae given in \cite{FIMU}, its lift to the nearly K\"ahler six-manifold is a Killing field for the sine-cone metric which does neither preserve $\Psi$ nor $\omega$ nor $J$.

\section{Left-invariant nearly $\ve$-K\"ahler structures on $\SL(2,\bR) \times \SL(2,\bR)$}
\label{mainsection}

\subsection{An algebraic prerequisite}
The following lemma is the key to proving the forthcoming structure result,
since it considerably reduces the number of algebraic equations on the nearly
$\ve$-K\"ahler candidates.

\bl \label{Cnormal}
Denote by $(\bR^{1,2}, \langle \cdot , \cdot \rangle)$  the vector space $\bR^3$
endowed with its standard Minkowskian scalar-product and 
denote by $\SO_0(1,2)$ the connected component of the identity of its group of isometries.  
Consider the action of $\SO_0(1,2) \times \SO_0(1,2)$ on the space of real  $3\times 3$ matrices $\mbox{\Mat}(3,\bR)$ given by
\bean
\Phi\,:\, \SO_0(1,2) \times \mbox{\Mat}(3,\bR) \times \SO_0(1,2) &\ra& \mbox{\Mat}(3,\bR)\\
(A,C,B) &\mapsto& A^t C B.
\eean
Then any invertible element $C \in \mbox{\Mat}(3,\bR)$ lies in the orbit of an element of the form
$$ 
\left(
\begin{array}{ccc}
 \alpha & x & y \\
0 & \beta & z \\
0 & 0& \gamma
\end{array} 
\right)
\quad 
\mbox{or}
\quad
\left(
\begin{array}{ccc}
0 & \beta & z \\
 \alpha & x & y \\
0 & 0& \gamma
\end{array} 
\right)
$$
with $ \alpha,\beta,\gamma,x,y,z \in \bR$ and $ \alpha \beta \gamma \ne 0.$
\el
\pf
Let an arbitrary invertible element $C \in  \mbox{\Mat}(3,\bR)$ be given. Denote by $\{e_1,e_2,e_3\}$ the standard basis of $\bR^{1,2}.$ There are three different cases:
\begin{itemize}
\item[1.)]
Suppose, that the first column $c$ of $C$ has negative length. We extend $c$ to a Lorentzian basis $\{l_1=c/\alpha,l_2,l_3\}$ with $\alpha:=\sqrt { |\langle c , c \rangle |}$. The linear map $L$ defined by extension of
$L(l_i)=e_i$ is by definition a Lorentz transformation. The transformation $L$
can be chosen time-oriented (by replacing $l_1$ by $\pm \,l_1$) and oriented 
(by replacing $l_3$ by $\pm\, l_3$). With this definition we obtain
$$\Phi(L^t,C,\id)= \left(
\begin{array}{cc}
\alpha & * \\
0 & C' \\
\end{array} 
\right) \mbox{ with an element } C' \in \mbox{\Mat}(2,\bR).
$$
Using the polar decomposition we can express 
$C'=O_1\,S$ as a product of $O_1 \in \SO(2)$ and a symmetric matrix $S$  in $\mbox{\Mat}(2,\bR)$ and diagonalise $S$ by $O_2 \in \SO(2).$ If we put
$$
L_1=\left(
\begin{array}{cc}
1 & 0 \\
0 & O_2^{-1} \, O_1^{-1}
\end{array} 
\right) 
\mbox{ and } L_2=\left(
\begin{array}{cc}
1&  0 \\
0 & O_2 
\end{array} \right)$$
 we obtain
$$\Phi(L_1^t,{\Phi}(L^t,C,\id),L_2)= \left(
\begin{array}{ccc}
\alpha & x & y \\
0 & \beta & 0 \\
0 & 0& \gamma
\end{array} 
\right).
$$
\item[2.)] 
Next suppose, that the first column $c$ of $C$ has positive length. Again, we extend $c$ to a Lorentzian basis 
$\{l_1,l_2=c/\alpha,l_3\}$ with $\alpha:=\sqrt { |\langle c , c \rangle |}$. The linear map $L$ defined by extension of $L(l_i)=e_i$ is by definition a Lorentz transformation. 
The transformation $L$ can be chosen  time-oriented (by replacing $l_1$ by
$\pm \,l_1$) and oriented (by replacing $l_3$ by $\pm\, l_3$). We get
$${\Phi}(L^t,C,\id)= \left(
\begin{array}{cc}
0 & * \\
\alpha & C' \\
0 & 
\end{array} 
\right) \mbox{ with an element } C' \in \mbox{\Mat}(2,\bR).
$$
The first column of this matrix is stable under the right-operation  of
$$
L_1=\left(
\begin{array}{cc}
1&  0 \\
0 & O_1 
\end{array} \right) \mbox{ with } O_1 \in \SO(2)$$
and there exists an element $O_1\in \SO(2)$ such
 that it holds
 $${\Phi}(\id,{\Phi}(L^t,C,\id),L_1)= \left(
\begin{array}{ccc}
0 & \beta & z \\
\alpha & x & y \\
0 & 0& \gamma
\end{array} 
\right).
$$
\item[3.)] Finally suppose, that it holds $\langle c , c \rangle =0.$ Then
  there exists an oriented and time-oriented Lorentz transformation $L$ 
such that $L(c)= \kappa (e_1 +e_2)$ with $\kappa \ne 0.$ Afterwards one finds
as in point 2.) an element $O\in \SO(2),$ such that it holds
$$C':=\Phi(L^t,C,O)= \left(
\begin{array}{ccc}
\kappa & c_1& *\\
\kappa & c_2 & *\\
0 & 0 & *
\end{array} 
\right). 
$$
Let
$$B(q) := \left(
\begin{array}{ccc}
\cosh(q) & \sinh(q)& 0\\
 \sinh(q) & \cosh(q) & 0\\
0 & 0 & 1
\end{array} 
\right).$$
{\bf Claim:} There exist $q_1,q_2 \in \bR$ such that 
$$ {\Phi}\left(\,B(q_1)^t\,,\,C'\,,\,B(q_2)\right)=\left(\begin{array}{ccc}
\alpha & x& y\\
0 & \beta & z\\
0 & 0 & \gamma
\end{array} 
\right). 
$$
To prove this claim let us first consider the right-action of $B(q)$ on $C'':= \Phi \left(B(q_1)^t, C', \id \right)$
$$ \Phi \left(\id, C'', B(q) \right)= \left(
\begin{array}{ccc}
c_{11}''\cosh(q) + c_{12}''\sinh(q) & *& *\\
c_{21}'' \cosh(q) + c_{22}''\sinh(q) &* & *\\
0 & 0 & *
\end{array} 
\right) \mbox{ for }  q \in \bR.$$
We choose $q_2$ such that $c_{21}'' \cosh(q_2) + c_{22}''\sinh(q_2)$ vanishes.
This is only possible if $-c_{22}''/c_{21}''$ is in the range of $\coth,$ i.e. $|c_{22}''/c_{21}''|> 1.$ \\
In the sequel we show, that this can always be achieved by the left-action of an element $B(q_1)$ on $C'$ and that $c_{21}'' \ne 0$. In fact, it is
\bean
c_{22}''&=&  c_1\sinh(q_1) +c_2\cosh(q_1)  \\
c_{21}''&=& \kappa(\sinh(q_1) + \cosh(q_1))= \kappa e^{q_1}\\
\frac{c_{22}''}{c_{21}''}&=& \frac{c_1+c_2}{2\kappa} + \frac{c_2-c_1}{2\kappa}e^{-2q_1}.
\eean
We observe, that $c_1\ne c_2,$ since the matrix $C$ is invertible. Therefore
we can always achieve $|{c_{22}''}/{c_{21}''}|>1.$ This proves the claim and
finishes the proof of the lemma. 
\end{itemize}
\end{proof}

\subsection{Proof of the uniqueness result}
Finally, we prove our main result which is the following theorem. By a homothety, we define the rescaling of the metric by a real number which we do not demand to be positive since we are working with all possible signatures.
\bt \label{main}
Let $G$ be a Lie group with Lie algebra $\mathfrak{sl}(2,\bR)$. 
Up to homothety, there is a unique left-invariant nearly $\ve$-K\"ahler structure with $\| \n \pJ \|^2 \ne 0$ on $G \times G$. This is the nearly pseudo-K\"ahler structure of signature (4,2) constructed as 3-symmetric space in the introduction. In particular, there is no left-invariant nearly para-K\"ahler structure.
\et
\br The proof also shows that there there is a left-invariant
nearly $\ve$-K\"ahler structure of non-zero type on $G \times H$ with $Lie(G) = Lie (H) = \mathfrak{sl}(2,\bR)$ if $G \ne H$ which is unique up to homothety \emph{and} exchanging the orientation. 
\er
\pf
More precisely, we will prove uniqueness up to equivalence of left-invariant almost $\ve$-Hermitian structures and homothety. We will consider the algebraic exterior system 
\begin{eqnarray}
d\o &=& 3 \, \psi^+,  \label{NKI} \\
d\psi^- &=& 2 \, \o \wedge \o \label{NKII}
\end{eqnarray}
on the Lie algebra $\mathfrak{sl}(2,\bR) \oplus \mathfrak{sl}(2,\bR)$. By Theorem \ref{ext_diff_syst_th}, solutions of this system are in one-to-one correspondence to left-invariant nearly $\ve$-K\"ahler structures on $G \times G$ with $\| \n \pJ \|^2 = 4$. This normalisation can always be achieved by applying a homothety. Furthermore, two solutions which are isomorphic under an inner Lie algebra automorphism from 
$$ {\rm Inn}(\mathfrak{sl}(2,\bR) \oplus \mathfrak{sl}(2,\bR))= \SO_0(1,2) \times \SO_0(1,2)$$ are equivalent under the corresponding Lie group isomorphism. Since both factors are equal, we can also lift the outer Lie algebra automorphism exchanging the two summands to the group level. In summary, it suffices to show the existence of a solution of the algebraic exterior system \eqref{NKI}, \eqref{NKII} on the Lie algebra which is unique up to inner Lie algebra automorphisms and exchanging the summands.

A further significant simplification is the observation that all tensors defining a nearly $\ve$-K\"ahler structure of non-zero type can be constructed out of the fundamental two-form $\omega$ with the help of the first nearly K\"ahler equation \eqref{NKI} and the stable form formalism described in section \ref{struct_red_ect}. We break the main part of the proof into three lemmas, step by step  simplifying $\omega$ under Lie algebra automorphisms in a fixed Lie bracket. 

We call $\{ e_1, e_2, e_3 \}$ a \emph{standard basis} of $\mathfrak{so}(1,2)$ if the Lie bracket satisfies
\begin{eqnarray*}
de^1= - e^{23}\, , \quad de^2= e^{31}\, , \quad de^3= e^{12}.
\end{eqnarray*}
In this basis, an inner automorphism in $\SO_0(1,2)$ acts by usual matrix multiplication on $\mathfrak{so}(1,2)$. 
\begin{lemma}
\label{lemmablockform} Let $\g = \h = \mathfrak{so}(1,2)$ and let $\omega$ be a non-degenerate two-form in $$\Lambda^2(\g \oplus \h)^* = \Lambda^2 \g^* \oplus (\g \otimes \h) \oplus \Lambda^2 \h^*.$$ Then we have 
\begin{eqnarray}
d\omega^2 = 0 \quad \Leftrightarrow \quad \omega \in \g \otimes \h.
\end{eqnarray}
\end{lemma}
\begin{proof}
By inspecting the standard basis, we observe that all two-forms on $\mathfrak{so}(1,2)$ are closed whereas no non-trivial 1-form is closed. Thus, when separately taking the exterior derivative of the components of $\omega^2$ in $\Lambda^4 = (\Lambda^3 \g^* \otimes \h^*) \oplus (\Lambda^2 \g^* \otimes \Lambda^2 \h^*) \oplus (\g^* \otimes \Lambda^3 \h^*$), the equivalence is easily deduced.
\end{proof}
\begin{lemma}\label{lemmaomeganormal}
Let $\g = \h = \mathfrak{so}(1,2)$ and let $\{ e^1, e^2, e^3 \}$ be a basis of $\g^*$ and $\{ e^4, e^5, e^6 \}$ a basis of $\h^*$ such that the Lie brackets are given by
\begin{equation}
\label{bas}
de^1= - e^{23} , \quad de^2= e^{31} , \quad de^3= \tau e^{12} \quad \mbox{and} \quad
de^4= -e^{56} , \quad de^5= e^{64} , \quad de^6= e^{45}
\end{equation}
for some $\tau \in \{ \pm 1\}$. Then, every non-degenerate two-form $\omega$  on $\g \oplus \h$ satisfying $d\omega^2 =0$ can be written 
\begin{eqnarray}
\label{omeganormal} \omega &=&\alpha \:e^{14}+\beta \:e^{25}+\gamma \:e^{36}+x \:e^{15}+y \:e^{16}+z \:e^{26}
\end{eqnarray}
for $\alpha, \beta, \gamma \in \bR - \{ 0 \}$ and $x,y,z \in \bR$ modulo an automorphism in $\SO_0(1,2) \times \SO_0(1,2)$.
\end{lemma}
\begin{proof}
We choose standard bases $\{ e^1, e^2, e^3 \}$ for $\g$ and $\{ e^4, e^5, e^6 \}$ for $\h$. Using the previous lemma and the assumption $d\omega^2 =0$, we may write $\omega = \sum_{i,j=1}^3 c_{ij} e^{i(j+3)}$ for an invertible matrix $C = \left(c_{ij}\right) \in \mbox{\Mat}(3,\bR)$. When a pair $(A,B) \in \SO_0(1,2) \times \SO_0(1,2)$ acts on the two-form $\omega$, the matrix $C$ is transformed to $A^t C B$. Applying Lemma \ref{Cnormal}, we can achieve by an inner automorphism that $C$ is in one of the normal forms given in that lemma. However, an exchange of the base vectors $e_1$ and $e_2$ corresponds exactly to exchanging the first and the second row of $C$. Therefore, we can always write $\omega$ in the claimed normal form by adding the sign $\tau$ in the Lie bracket of the first summand $\g$.
\end{proof}

\begin{lemma}
Let $\{ e^1, \dots, e^6 \}$ be a basis of $\mathfrak{so}(1,2) \times \mathfrak{so}(1,2)$ such that 
\begin{equation}
\label{bas2}
de^1= - e^{23}\, , \quad de^2= e^{31}\, , \quad de^3= e^{12} \quad \mbox{and} \quad
de^4= - e^{56}\, , \quad de^5= e^{64}\, , \quad de^6= e^{45}. 
\end{equation}
Then the only $\SU^\ve(p,q)$-structure $(\omega,\psi^+)$ modulo inner automorphisms and modulo exchanging the summands, which solves the two nearly $\ve$-K\"ahler equations \eqref{NKI} and \eqref{NKII}, is determined by 
\begin{equation}
  \label{omegaNK}
   \omega = \frac{\sqrt 3}{18} (e^{14} + e^{25} + e^{36}). 
\end{equation}
\end{lemma}
\begin{proof}
Since $d\omega^2=0$ by the second equation \eqref{NKII}, we can choose a basis satisfying \eqref{bas} such that $\omega$ is in the normal form \eqref{omeganormal}. In order to satisfy the first equation \eqref{NKI}, we have to set 
\begin{eqnarray*}
 3 \psi^+ = d\omega &=& -\alpha \:e^{234} +\alpha \:e^{156} -x \:e^{235} +x \:e^{146} -y \:e^{236} -y \:e^{145}\\ &-&\beta \:e^{135} +\beta \:e^{246} -z \:e^{136} -z \:e^{245} +\tau \gamma \:e^{126} -\gamma \:e^{345}.
\end{eqnarray*}
The compatibility $\omega \wedge \psi^+ = 0$ is equivalent to $d(\omega^2)=0.$ It remains to determine all solutions of the second nearly $\ve$-K\"ahler equation \eqref{NKII} modulo automorphisms.

For the sake of readability, we identify $\Lambda^6 (\g \oplus \h)^*$ with $\bR$ by means of $e^{123456}$. Supported by Maple, we compute
\begin{eqnarray*}
K_{\psi^+}(e_1)&=&
(x^2+y^2+z^2-\alpha^2 +\beta^2+\tau \gamma^2) e_1
-(2 x \beta+ 2 y z) e_2\\
&-&2 \tau \gamma y e_3
+2 \tau \gamma  \beta e_4,\\
K_{\psi^+}(e_2)&=&
(2 x \beta+2 y z) e_1
+(-x^2-y^2-z^2+\alpha^2-\beta^2+\tau \gamma^2) e_2\\
&-&2 \tau \gamma z e_3
+2 \tau \gamma x e_4
-2 \tau \alpha \gamma e_5,\\
K_{\psi^+}(e_3)&=&
2 y \gamma e_1
-2 z \gamma e_2
+(-x^2-y^2+z^2+\alpha^2+\beta^2-\tau \gamma^2) e_3\\
&+&(2 y \beta -2 x z) e_4
+2 \alpha z e_5
-2 \alpha \beta e_6,\\
K_{\psi^+}(e_4)&=&
-2 \beta \gamma e_1
+2 x \gamma e_2
+(2 y \beta-2 x z) e_3\\
&+&(x^2+y^2-z^2+\alpha^2-\beta^2-\tau \gamma^2) e_4
-2 \alpha x e_5
-2 \alpha y e_6,\\
K_{\psi^+}(e_5)&=&
2 \alpha \gamma e_2
-2 \alpha z e_3
+2 \alpha x e_4\\
&+&(-x^2+y^2-z^2-\alpha^2+\beta^2-\tau \gamma^2) e_5
+(2 \beta z -2 x y) e_6,\\
K_{\psi^+}(e_6)&=&
2 \alpha \beta e_3
+2 \alpha y e_4\\
&+&(2 \beta z -2 x y) e_5
+(x^2-y^2+z^2-\alpha^2-\beta^2+\tau \gamma^2) e_6.
\end{eqnarray*}
We assume that $\lambda(\psi^+) \ne 0$ and check this a posteriori for the solutions we find. Hence, we can set $k:= \frac{1}{\pm \sqrt{|\lambda(\psi^+)|}}$ and $J_{\psi^+} = k K_{\psi^+} $. Since $\psi^+ + \hati J^*_{\psi^+} \psi^+$ is a $(3,0)$-form with respect to $J_{\psi^+}$, we have $\psi^- = J^*_{\psi^+} \psi^+ = \varepsilon \psi^+ (J_{\psi^+} . , ., .)$ which turns out to be
\begin{eqnarray*}
\varepsilon\frac{27}{k}\, \psi^- &=& 
2 \tau \alpha \beta \gamma \:e^{123} 
+ 2 \tau y \alpha \gamma \:e^{124} 
+ 2 \tau \gamma (x y-\beta z) \:e^{125} 
- 2 (x \beta+y z) \alpha \:e^{134} \\
&+& \tau  \gamma (- x^2 +  y^2 -  z^2 +  \alpha^2 +  \beta^2 - \tau \gamma^2) \:e^{126}\\
&-& \{ \beta ( x^2- y^2 - z^2 + \alpha^2 -\beta^2 + \tau \gamma^2) + 2 xyz \}\:e^{135} \\
&+&\{ z (x^2 - y^2 + z^2 - \alpha^2 + \beta^2 + \tau \gamma^2) - 2 x y \beta \} \:e^{136} \\
&-&\{ y (-x^2-y^2+ z^2 + \alpha^2-\beta^2+\tau \gamma^2) + 2 x z \beta \} \:e^{145} \\
&-&\{ x (x^2+ y^2+ z^2 - \alpha^2-\beta^2+\tau \gamma^2) - 2 y z \beta \} \:e^{146} \\
&-&\alpha (x^2+ y^2+ z^2 -\alpha^2+ \beta^2 + \tau \gamma^2) \:e^{156} \\
&-&\alpha (x^2+ y^2+ z^2 -\alpha^2+ \beta^2 + \tau \gamma^2) \:e^{234} \\
&-&\{ x(x^2+ y^2+ z^2 - \alpha^2 -\beta^2 +\tau \gamma^2) - 2 y z \beta \} \:e^{235} \\
&+&\{ y(-x^2-y^2+ z^2 + \alpha^2 -\beta^2 +\tau \gamma^2) + 2 x z \beta \} \:e^{236} \\
&-&\{ z(x^2 - y^2 + z^2 - \alpha^2 + \beta^2 + \tau \gamma^2) - 2 x y \beta \} \:e^{245} \\
&-&\{ \beta ( x^2- y^2 - z^2 + \alpha^2 -\beta^2 + \tau \gamma^2) + 2 x y z \} \:e^{246} \\
&+&\gamma (- x^2 +  y^2 -  z^2 +  \alpha^2 + \beta^2 - \tau \gamma^2) \:e^{345} \\
&-& 2 (x \beta+y z) \alpha \:e^{256}
- 2 \gamma (x y-\beta z) \:e^{346} 
- 2 y \alpha \gamma \:e^{356} +2 \alpha \beta \gamma \:e^{456}.
\end{eqnarray*}
Furthermore, we compute the exterior derivative 
\begin{eqnarray*}
  \varepsilon \frac{27}{k} d\psi^- &=&  
- 4 \tau \gamma \alpha y \:e^{1256} \:
- 4 \tau \gamma (x y - \beta z) \:e^{1246} \:
+ 4  \alpha (x \beta + y z) \:e^{1356} \\
&+&2 \tau \gamma (- x^2 +  y^2 -  z^2 +  \alpha^2 +  \beta^2 - \tau \gamma^2) \:e^{1245} \\
&+& 2 \{ \beta ( x^2- y^2 - z^2 + \alpha^2 -\beta^2 + \tau \gamma^2) + 2 x y z \} \:e^{1346} \\
&+& 2 \{ z (x^2 + y^2 +z^2 - \alpha^2+\beta^2+\tau \gamma^2 ) - 2 x y \beta \} \:e^{1345} \\
&+& 2 \{ y (- x^2-y^2+ z^2 + \alpha^2-\beta^2+\tau \gamma^2) + 2 x z \beta \} \:e^{2345} \\
&+& 2 \{ x(x^2+ y^2+ z^2 - \alpha^2- \beta^2 + \tau \gamma^2)-2 y z \beta \} \:e^{2346} \\
&+& 2 \alpha (x^2+ y^2+ z^2 -\alpha^2+ \beta^2 + \tau \gamma^2) \:e^{2356} 
\end{eqnarray*}
and
\begin{eqnarray*}
\omega^2 &=& 2 ( ( y \beta - xz) \: e^{1256} -  \alpha z \:e^{1246} - x \gamma \:e^{1356} -  \alpha \beta \:e^{1245} - \alpha \gamma \:e^{1346} -  \beta \gamma \:e^{2356}).
\end{eqnarray*}
The second nearly K\"ahler equation \eqref{NKII} is therefore equivalent to the following nine coefficient equations:
\begin{eqnarray*}
(\alpha \beta - 27 \ve k^{-1} \gamma ) \, x =& - \alpha \, y z, & \qquad \mbox{($e^{1356}$)} \\ 
(\tau \gamma \alpha - 27 \ve k^{-1} \beta)\, y =& - 27 \ve k^{-1} x z,  & \qquad \mbox{($e^{1256}$)} \\
(\tau \beta \gamma - 27 \ve k^{-1} \alpha)\, z =& \tau \gamma  x y,  & \qquad \mbox{($e^{1246}$)} \\ 
x^2+ y^2+ z^2 -\alpha^2+ \beta^2 + \tau \gamma^2 \, =& 54 \ve k^{-1} \frac{\beta \gamma}{\alpha},  & \qquad \mbox{($e^{2356}$)}  \\
z (x^2 + y^2 +z^2 - \alpha^2+\beta^2+\tau \gamma^2 ) =& - 2 \beta y x,  & \qquad \mbox{($e^{1345}$)} \\
 x^2- y^2 - z^2 + \alpha^2 -\beta^2 + \tau \gamma^2 \, =& 54 \ve k^{-1} \frac{\alpha \gamma}{\beta} - 2 \frac{x y z}{\beta},  & \qquad \mbox{($e^{1346}$)} \\ 
y (- x^2-y^2+ z^2 + \alpha^2-\beta^2+\tau \gamma^2) =& - 2 \beta z x,  & \qquad \mbox{($e^{2345}$)} \\ 
- x^2 +  y^2 -  z^2 +  \alpha^2 +  \beta^2 - \tau \gamma^2 \, =& 54 \tau \ve k^{-1} \frac{\alpha \beta}{\gamma},  & \qquad \mbox{($e^{1245}$)} \\
 x(-x^2- y^2- z^2 + \alpha^2+ \beta^2 - \tau \gamma^2) =& - 2 \beta y z.  & \qquad \mbox{($e^{2346}$)}
\end{eqnarray*}
Recall that $\alpha,\beta,\gamma \ne 0$ because $\o$ is non-degenerate.
We claim that there is no solution if any of $x$, $y$ or $z$ is different from zero.

On the one hand, assume that one of them is zero. Using one of the first three equations respectively, we find that at least one of the other two has to be zero as well. However, in all three cases, we may easily deduce that the third one has to be zero as well by comparing equations 4 and 5 respectively 6 and 7 respectively 8 and 9.

On the other hand, if we assume that all three of them are different from zero, the bracket in the first equation is necessarily different from zero and we may express $x$ by a multiple of $yz$. Substituting this expression into equations 2 and 3, yields expressions for $y^2$ and $z^2$ in terms of $\alpha$,$\beta$, $\gamma$ and $k$. But if we insert all this into equation 4 (or 6 or 8 alternatively), we end up with a contradiction after a slightly tedious calculation.

To conclude, we can set $x = y = z = 0$ without losing any solutions of the second nearly K\"ahler equation which simplifies to the equations
\begin{eqnarray*}
\alpha^3- \alpha \beta^2 - \tau  \alpha \gamma^2 - 54 \ve k^{-1} \beta \gamma &=& 0, \\
\beta^3 - \tau \beta \gamma^2 -  \beta \alpha^2 - 54 \ve k^{-1} \gamma \alpha &=& 0, \\
\gamma^3 - \tau \gamma \alpha^2 - \tau \gamma \beta^2 - 54 \ve k^{-1} \alpha \beta &=& 0 .
\end{eqnarray*}
Setting $c_1=\alpha^2 + \beta^2 + \tau \gamma^2$ and $c_2=54 \ve k^{-1} \alpha \beta \gamma$, these are equivalent to
\begin{eqnarray}
2 \alpha^4 - c_1 \alpha^2 - c_2 &=& 0, \nonumber\\
2 \beta^4 - c_1 \beta^2 - c_2 &=& 0, \label{sys}\\
2 \gamma^4 - c_1 \tau \gamma^2 - c_2 &=& 0 \nonumber.
\end{eqnarray}

To finish the proof, we have to show that all real solutions of the system \eqref{sys} are isomorphic under $\SO_0(1,2) \times \SO_0(1,2)$ to $$\alpha=\beta=\gamma= \frac{\sqrt 3}{18}\: , \qquad \tau=1.$$
Since $\alpha^2$, $\beta^2$ and $\tau \gamma^2$ satisfy the same quadratic
equation, at least two of them have to be identical, say
$\alpha^2=\beta^2$. However, if $\tau \gamma^2$ was the other root of the
quadratic equation, we would have $\alpha^2 + \tau \gamma^2 = \frac{1}{2} c_1$
and by definition of $c_1$ at the same time $2 \alpha^2 + \tau \gamma^2 =
c_1$. This would only be possible, if $\gamma$ was zero, a contradiction to
the non-degeneracy of $\omega$. Therefore $\tau$ has to be $+1$ and
$\alpha$,$\beta$ and $\gamma$ have to be identical up to sign. 
By applying one of the proper and orthochronous Lorentz transformations 
\begin{eqnarray*}
  \begin{pmatrix}
    1&0&0\\0&-1&0\\0&0&-1
  \end{pmatrix}, \qquad 
  \begin{pmatrix}
    1&0&0\\0&0&1\\0&-1&0
  \end{pmatrix}, \qquad
  \begin{pmatrix}
    1&0&0\\0&0&-1\\0&1&0
  \end{pmatrix}
\end{eqnarray*}
on, say, the second summand, it is always possible to achieve that the signs of $\alpha$, $\beta$ and $\gamma$ are identical. 

So far, we found a basis satisfying \eqref{bas2} such that 
$$ \omega = \alpha (e^{14} + e^{25} + e^{36}). $$
It is straightforward to check that the quartic invariant in this basis is
\begin{equation}
  \label{lambdanegativ}
\lambda ( \frac{1}{3} d \omega) = -\frac{1}{27}\alpha^4.  
\end{equation}
Therefore, there cannot exist a nearly para-K\"ahler structure and we can set $\ve = -1$. Inserting $k= \pm \frac{1}{\sqrt{-\lambda}} = \pm 3 \sqrt 3 \alpha^{-2}$ into equations \eqref{sys} yields
$$ 2 \alpha^4 - 3 \alpha^4 \pm \frac{54}{3 \sqrt 3} \alpha^5 = 0 \iff \alpha = \pm \frac{1}{18} \sqrt{3}.$$
Finally, we can achieve that $\alpha$ is positive by applying the Lie algebra automorphism exchanging the two summands, i.e. $e_i \mapsto e_{i+3 \; {\rm mod} \; 6}$ and the lemma is proven.
\end{proof}
In fact, the uniqueness, existence and non-existence statements claimed in the theorem follow directly from this lemma and formula \eqref{lambdanegativ}. 

As explained in the introduction, we know that there is a left-invariant nearly pseudo-K\"ahler structure of indefinite signature on all the groups in question. After applying a homothety, we can achieve $\| \n \pJ \|^2 = 4$ and this structure has to coincide with the unique structure we just constructed. Therefore, the indefinite metric has to be of signature (4,2) by our sign conventions. 

We summarise the data of the unique nearly pseudo-K\"ahler structure in the basis \eqref{bas2} and can easily double-check the signature of the metric explicitly:
\begin{eqnarray*}
\omega &=& \frac{1}{18} \sqrt{3}\; (e^{14}+ \:e^{25}+ \:e^{36})\\
\psi^+ &=& \frac{1}{54} \sqrt{3} \; (e^{126} -\:e^{135} +\:e^{156} -\:e^{234} +\:e^{246}  -\:e^{345})\\
\psi^- &=& - \frac{1}{54} (2\:e^{123} + \:e^{126} - \:e^{135} - \:e^{156} - \:e^{234} - \:e^{246} + \:e^{345} + 2\:e^{456} )\\
J(e_1) &=& - \frac{1}{3} \sqrt 3 \; e_1 - \frac{2}{3} \sqrt 3 \; e_4 \, , \quad 
J(e_4) =\:\:\:\, \frac{2}{3} \sqrt 3 \; e_1 + \frac{1}{3} \sqrt 3 \; e_4\\
J(e_2) &=& - \frac{1}{3} \sqrt 3 \; e_2 + \frac{2}{3} \sqrt 3 \; e_5 \, , \quad 
J(e_5) = - \frac{2}{3} \sqrt 3 \; e_2 + \frac{1}{3} \sqrt 3 \; e_5\\
J(e_3) &=& - \frac{1}{3} \sqrt 3 \; e_3 + \frac{2}{3} \sqrt 3 \; e_6 \, , \quad 
J(e_6) = - \frac{2}{3} \sqrt 3 \; e_3 + \frac{1}{3} \sqrt 3 \; e_6\\
g &=& \frac{1}{9} \; ( \; (e^1)^2 - \; (e^2)^2 - \; (e^3)^2 + \; (e^4)^2 - \; (e^5)^2 - \; (e^6)^2 - \; e^1 \cdot e^4 - \; e^2 \cdot e^5 - \; e^3 \cdot e^6).
\end{eqnarray*}
\epf

Observing that in \cite{Bu} very similar arguments have been applied to the Lie group $S^3 \times S^3$, we find the following non-existence result. 

\bp
On the Lie groups $G \times H$ with $Lie (G) = Lie(H) = \mathfrak{so}(3)$, there is neither a left-invariant nearly para-K\"ahler structure of non-zero type nor a left-invariant nearly pseudo-K\"ahler structure with an indefinite metric.
\ep
\pf The unicity of the left-invariant nearly K\"ahler structure $S^3 \times S^3$ is proved in \cite{Bu}, section 3, with a strategy analogous to the proof of Theorem \ref{main}. In the following, we will refer to the 
English version \cite{Bu2}. There, it is shown in the proof of Proposition 2.5, that for any solution of the exterior system 
\begin{eqnarray*}
  d \o &=& 3 \psi^+ \\
  d \psi^+ &=& -2 \mu \o ^2
\end{eqnarray*}
 there is a basis of the Lie algebra of $S^3 \times S^3$ and a real constant $\alpha$ such that 
\begin{eqnarray*}
&de^1= e^{23}\, , \quad de^2= e^{31}\, , \quad de^3= e^{12} \quad \mbox{and} \quad
de^4= e^{56}\, , \quad de^5= e^{64}\, , \quad de^6= e^{45},& \\
&\omega = \alpha (e^{14} + e^{25} + e^{36}).&
\end{eqnarray*}
In this basis, a direct computation or formula (18) in \cite{Bu2} show that the quartic invariant that we denote by $\lambda$ is 
\begin{eqnarray*}
\lambda = - \frac{1}{27} \alpha^4
\end{eqnarray*}
with respect to the volume form $e^{123456}$. Therefore, a nearly para-K\"ahler structure cannot exist on all the Lie groups with the same Lie algebra as $S^3 \times S^3$ by Theorem \ref{ext_diff_syst_th}. A nearly pseudo-K\"ahler structure with an indefinite metric cannot exist either, since the induced metric is always definite as computed in the second part of Lemma 2.3 in \cite{Bu2}.
\epf

\end{document}